\documentclass[10pt]{amsart}
\input{epsf}
\usepackage{amssymb,latexsym}
\topskip  \newcommand{\resp}[1]{{\rm [}res\-pec\-ti\-ve\-ly; #1{\rm ]}}

\numberwithin{equation}{section}
\newtheorem{theorem}{Theorem}[section]
\newtheorem{proposition}[theorem]{Proposition}

\newtheorem{conjecture}[theorem]{Conjecture}

\newtheorem{lemma}[theorem]{Lemma}

\begin{document}

\pagenumbering{arabic}
\pagestyle{headings}
\def\sof{\hfill\rule{2mm}{2mm}}
\def\ls{\leq}
\def\gs{\geq}
\def\SS{\frak S}
\def\qq{{\bold q}}
\def\txx{{\frac1{2\sqrt{x}}}}

\title{Counting occurrences of $132$ in an even permutation}
\maketitle

\begin{center}Toufik Mansour \footnote{Research financed by EC's
IHRP Programme, within the Research Training Network "Algebraic
Combinatorics in Europe", grant HPRN-CT-2001-00272}
\end{center}

\begin{center}{Department of Mathematics, Chalmers University of
Technology, S-41296 G\"oteborg, Sweden

        {\tt toufik@math.chalmers.se} }
\end{center}

\def\EE1{-\frac{1}{2}(1-2x-x^2)+\frac{1-3x}{4\sqrt{1-4x}}+\frac{1-3x-4x^2+4x^3}{4\sqrt{1-4x^2}}}
\def\OO1{-\frac{x}{2}(1+x)+\frac{1-3x}{4\sqrt{1-4x}}-\frac{1-3x-4x^2+4x^3}{4\sqrt{1-4x^2}}}
\def\sign{{\rm sign}}
\def\vv{{\bf v}}
\def\rr{{\bf r}}
\def\uu{{\bf u}}
\def\vr{\varnothing}
\def\BB{\mathcal B}

\section*{Abstract}
We study the generating function for the number of even (or odd)
permutations on $n$ letters containing exactly $r\gs0$ occurrences
of $132$. It is shown that finding this function for a given $r$
amounts to a routine check of all permutations in $S_{2r}$.

\noindent {2000 Mathematics Subject Classification}: Primary
05A05, 05A15; Secondary 05C90
\section{Introduction}
Let $[n]=\{1,2,\ldots,n\}$ and $\SS_n$ denote the set of all
permutations of $[n]$. We shall view permutations in $\SS_n$ as
words with $n$ distinct letters in $[n]$. A pattern is a
permutation $\sigma\in\SS_k$, and an occurrence of $\sigma$ in a
permutation $\pi=\pi_1\pi_2 \cdots\pi_n\in\SS_n$ is a subsequence
of $\pi$ that is order equivalent to $\sigma$. For example, an
occurrence of $132$ is a subsequence $\pi_i\pi_j\pi_k$ ($1\leq i <
j < k\leq n$) of $\pi$ such that $\pi_i<\pi_k<\pi_j$. We denote by
$\tau(\pi)$ the number of occurrences of $\tau$ in $\pi$, and we
denote by $s_{\sigma}^r(n)$ the number of permutations
$\pi\in\SS_n$ such that $\sigma(\pi)=r$.

In the last decade much attention has been paid to the problem of
finding the numbers $s_{\sigma}^r(n)$ for a fixed $r\geq 0$ and a
given pattern $\tau$ (see
\cite{AlFr00,At99,Bo97b,Bo97a,Bo98,ChWe99,Ma00,MaVa01,NoZe96,
Ro99,SiSc85,St94,St96,We95}). Most of the authors consider only
the case $r=0$, thus studying permutations {\it avoiding\/} a
given pattern. Only a few papers consider the case $r>0$, usually
restricting themselves to patterns of length $3$.  Using two
simple involutions (\emph{reverse} and \emph{complement}) on
$\SS_n$ it is immediate that with respect to being
equidistributed, the six patterns of length three fall into the
two classes $\{123,321\}$ and $\{132,213,231,312\}$. Noonan
\cite{No96} proved that $$s_{123}^1(n)=\frac 3n\binom{2n}{n-3}.$$
A general approach to the problem was suggested by Noonan and
Zeilberger \cite{NoZe96}; they gave another proof of Noonan's
result, and conjectured that
        $$s_{123}^2(n)=\frac{59n^2+117n+100}{2n(2n-1)(n+5)}\binom{2n}{n-4}$$
and
$$s_{132}^1(n)=\binom{2n-3}{n-3}.$$
The first conjecture was proved by Fulmek \cite{Fulmek} and the
latter conjecture was proved by B\'ona in \cite{Bo98}. A
conjecture of Noonan and Zeilberger states that $s_{\sigma}^r(n)$
is $P$-recursive in $n$ for any $r$ and $\tau$. It was proved by
B\'ona \cite{Bo97c} for $\sigma=132$. Mansour and Vainshtein
\cite{MaVa01} suggested a new approach to this problem in the case
$\sigma=132$, which allows one to get an explicit expression for
$s_{132}^r(n)$ for any given $r$. More precisely, they presented
an algorithm that computes the generating function $\sum_{n\geq0}
s_{132}^r(n)x^n$ for any $r\geq0$. To get the result for a given
$r$, the algorithm performs certain routine checks for each
element of the symmetric group $\SS_{2r}$. The algorithm has been
implemented in C, and yields explicit results for $1\leq r\leq 6$.

Let $\pi$ be any permutation. The number of {\em inversions} of
$\pi$ is given by $i_\pi=|\{(i,j): \pi_i>\pi_j,\ i<j\}|$. The {\em
signature} of $\pi$  is given by $\sign(\pi)=(-1)^{i_\pi}$. We say
$\pi$ is an {\em even permutation} \resp{{\em odd permutation}} if
$\sign(\pi)=1$ \resp{$\sign(\pi)=-1$}. We denote by $E_n$
\resp{$O_n$} the set of all even \resp{odd} permutations in
$\SS_n$. Clearly, $|E_n|=|O_n|=\frac{1}{2}n!$ for all $n\geq2$.
The following lemma holds immediately by definitions.
\begin{lemma}\label{slem1}
For any permutation $\pi$,  $\sign(\pi)=(-1)^{21(\pi)}$.
\end{lemma}

We denote by $e_{\sigma}^r(n)$ \resp{$o_{\sigma}^r(n)$} the number
of even \resp{odd} permutations $\pi\in E_n$ \resp{$\pi\in O_n$}
such that $\sigma(\pi)=r$.

Apparently, for the first time the relation between even (odd)
permutations and pattern avoidance problem was suggested by Simion
and Schmidet in \cite{SiSc85} for $\sigma\in\SS_3$. In
particularly, Simion and Schmidt \cite{SiSc85} proved that
$$e_{132}^0(n)=\frac{1}{2(n+1)}\binom{2n}{n}+\frac{1}{n+1}\binom{n-1}{(n-1)/2}\mbox{
and }
o_{132}^0(n)=\frac{1}{2(n+1)}\binom{2n}{n}-\frac{1}{n+1}\binom{n-1}{(n-1)/2}.$$

In this paper, as a consequence of \cite{MaVa02}, we suggest a new
approach to this problem in the case of even (or odd) permutations
where $\sigma=132$, which allows one to get an explicit expression
for $e_{132}^r(n)$ for any given $r$. More precisely, we present
an algorithm that computes the generating functions
$E_r(x)=\sum_{n\geq0} e_{132}^r(n)x^n$ and $O_r(x)=\sum_{n\geq0}
o_{132}^r(n)x^n$ for any $r\geq0$. To get the result for a given
$r$, the algorithm performs certain routine checks for each
element of the symmetric group $\SS_{2r}$. The algorithm has been
implemented in C, and yields explicit results for $0\leq r\leq 6$.
\section{Recall definitions and preliminary results}
\begin{center}
\begin{figure}[h]
\epsfxsize=2.5in \epsffile{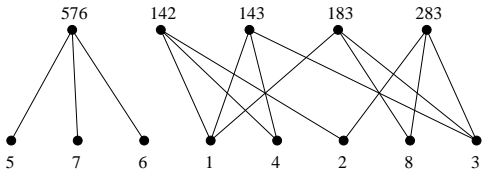} \caption{Graph $G_\pi$ for
$\pi=57614283$} \label{graph}
\end{figure}
\end{center}

To any $\pi\in\SS_n$ we assign a bipartite graph $G_\pi$ in the
following way. The vertices in one part of $G_\pi$, denoted $V_1$,
are the entries of $\pi$, and the vertices of the second part,
denoted $V_3$, are the occurrences of $132$ in $\pi$. Entry $i\in
V_1$ is connected by an edge to occurrence $j\in V_3$ if $i$
enters $j$. For example, let $\pi=57614283$, then $\pi$ contains
$5$ occurrences of $132$, and the graph $G_\pi$ is presented on
Figure~\ref{graph}.
We denote by $G_\pi^n$ the connected component of $G_\pi$
containing entry $n$. Let $\pi(i_1),\dots,\pi(i_s)$ be the entries
of $\pi$ belonging to $G_\pi^n$, and let
$\sigma=\sigma_\pi\in\SS_s$ be the corresponding permutation. We
say that $\pi(i_1),\dots,\pi(i_s)$ is the {\it kernel\/} of $\pi$
and denote it $\ker\pi$, $\sigma$ is called the {\it shape\/} of
the kernel, or the {\it kernel shape}, $s$ is called the {\it
size\/} of the kernel, and the number of occurrences of $132$ in
$\ker\pi$ is called the {\it capacity\/} of the kernel. For
example, for $\pi=57614283$ as above, the kernel equals $14283$,
its shape is $14253$, the size equals $5$, and the capacity equals
$4$.

\begin{theorem}\label{MVth1} {\rm (Mansour and
Vainshtein~\cite[Theorem~1]{MaVa02})} Let $\pi\in\SS_n$ such that
$132(\pi)=r$, then the size of the kernel of $\pi$ is at most
$2r+1$.
\end{theorem}

\begin{center}
\begin{figure}[h]
\epsfxsize=2.5in \epsffile{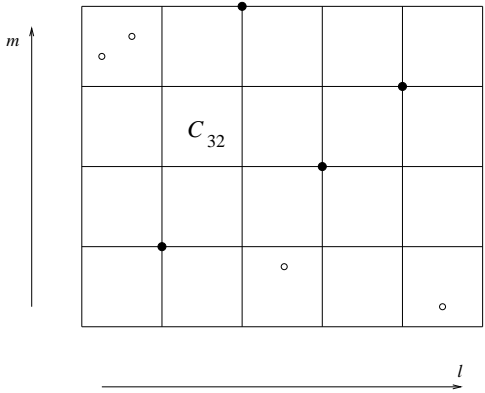} \caption{Kernel cell
decomposition for $\pi\in\SS(1423)$} \label{cell}
\end{figure}
\end{center}

We say that $\rho$ is a {\it kernel permutation\/} if it is the
kernel shape for some permutation $\pi$. Evidently $\rho$ is a
kernel permutation if and only if $\sigma_\rho=\rho$. Let
$\rho\in\SS_s$ be an arbitrary kernel permutation. We denote by
$\SS(\rho)$ the set of all the permutations of all possible sizes
whose kernel shape equals $\rho$. For any $\pi\in\SS(\rho)$ we
define the {\it kernel cell decomposition\/} as follows. The
number of cells in the decomposition equals $s(s+1)$. Let
$\ker\pi=\pi(i_1),\dots,\pi(i_s)$; the {\it cell\/}
$C_{ml}=C_{ml}(\pi)$ for $1\ls l\ls s+1$ and $1\ls m\ls s$ is
defined by
$$
C_{ml}(\pi)=\{\pi(j)\: i_{l-1}<j<i_{l},
\;\pi(i_{\rho^{-1}(m-1)})<\pi(j)< \pi(i_{\rho^{-1}(m)})\},
$$
where $i_0=0$, $i_{s+1}=n+1$,  and $\pi(0)=0$ for any $\pi$. If
$\pi$ coincides with $\rho$ itself, then all the cells in the
decomposition are empty. An arbitrary permutation in $\SS(\rho)$
is obtained  by filling in some of the cells in the cell
decomposition. A cell $C$ is called {\it infeasible\/} if the
existence of an entry $a\in C$ would imply an occurrence of $132$
that contains $a$ and two other entries $x,y\in\ker\pi$. Clearly,
all infeasible cells are empty for any $\pi\in\SS(\rho)$. All the
remaining cells are called {\it feasible\/}; a feasible cell may,
or may not, be empty. Consider the permutation $\pi=67382451$. The
kernel of $\pi$ equals $3845$, its shape is $1423$. The cell
decomposition of $\pi$ contains four feasible cells:
$C_{13}=\{2\}$, $C_{14}=\varnothing$, $C_{15}=\{1\}$, and
$C_{41}=\{6,7\}$ (see Figure~\ref{cell}). All the other cells are
infeasible; for example, $C_{32}$ is infeasible, since if $a\in
C_{32}$, then $a\pi'(i_2)\pi'(i_4)$ is an occurrence of $132$ for
any $\pi'$ whose kernel is of shape $1423$.

As another example, permutation $\widetilde
\pi=11\,10\,7\,12\,4\,6\,5\,8\,3\,9\,1\,2$ belongs to the same
class $\SS(1423)$. Its kernel is $7\,12\,8\,9$, and the feasible
cells are $C_{13}=\{4,6,5\}$, $C_{14}=\{3\}$, $C_{15}=\{1,2\}$,
$C_{41}=\{11,10\}$.

Given a cell $C_{ij}$ in the kernel cell decomposition, all the
kernel entries can be positioned with respect to $C_{ij}$. We say
that $x=\pi(i_k)\in\ker\pi$ lies {\it below\/} $C_{ij}$ if
$\rho(k)<i$, and {\it above\/} $C_{ij}$ if $\rho(k)\gs i$.
Similarly, $x$ lies to the {\it left\/} of $C_{ij}$ if $k<j$, and
to the {\it right\/} of $C_{ij}$ if $k \gs j$. As usual, we say
that $x$ lies to the {\it southwest\/} of $C_{ij}$ if it lies
below $C_{ij}$ and to the left of it; the other three directions,
northwest, southeast, and northeast, are defined similarly. Let us
define a partial order $\prec$ on the set of all feasible cells by
saying that $C_{ml}\prec C_{m'l'}\ne C_{ml}$ if $m\gs m'$ and
$l\ls l'$.

\begin{theorem}\label{MVlem2}
{\rm (Mansour and Vainshtein~\cite{MaVa02})}

{\rm (i)} $\prec$  is a linear order.

{\rm(ii)} Let $C_{ml}$ and $C_{ml'}$ be two nonempty feasibly
cells such that $l<l'$. Then for any pair of entries $a\in
C_{ml}$, $b\in C_{ml'}$, one has $a>b$.

{\rm (iii)} Let $C_{ml}$ and $C_{m'l}$ be two nonempty feasibly
cells such that $m<m'$. Then any entry $a\in C_{ml}$ lies to the
right of any entry $b\in C_{m'l}$.
\end{theorem}
%
%
Let $\pi$ be any permutation with a kernel permutation $\rho$, and
assume that the feasible cells of the kernel cell decomposition
associated with $\rho$ are ordered linearly according to $\prec$,
$C^1,C^2,\ldots,C^{f(\rho)}$. Let $d_j$ be the size of $C^j$. For
example, let $\pi=67382451$ with kernel permutation $\rho=1423$,
as on Figure~\ref{cell}, then $d_1=2$, $d_2=1$, $d_3=0$, and
$d_4=1$.\\
We denote by $l_j(\rho)$ the number of the entries of $\rho$ that
lie to the north-west from $C^j$ or lie to the south-east from
$C^j$. For example, let $\rho=1423$, as on Figure~\ref{cell}, then
$l_1(\rho)=3$, $l_2(\rho)=2$, $l_3(\rho)=3$, and $l_4(\rho)=4$.
Clearly, $l_1(\rho)=s(\rho)-1$ and $l_{f(\rho)}=s(\rho)$ for any
nonempty kernel permutation $\rho$.

\begin{lemma}\label{alem1}
For any permutation $\pi$ with a kernel permutation $\rho$,
$$\sign(\pi)=(-1)^{\left(\sum_{1\leq i\leq j\leq
f(\rho)}\;d_id_j+\sum_{j=1}^{f(\rho)}\;d_jl_j(\rho)\right)}\cdot
\sign(\rho)\cdot\prod_{j=1}^{f(\rho)}\sign(C^j).$$
\end{lemma}
\begin{proof}
To verify this formula, let us count the number of occurrences of
the pattern $21$ in $\pi$. There four possibilities for an
occurrence of $21$ in $\pi$. The first possibility is an
occurrence occurs in one of the cells $C^j$, so in this case there
are $\sum_{j=1}^{f(\rho)}21(C^j)$ occurrences. The second
possibility is an occurrence occurs in $\rho$, so there are
$21(\rho)$ occurrences. The third possibility is an occurrence of
two elements which the first belongs to $\rho$ and the second
belongs to $C^i$, so there are $\sum_{j=1}^{f(\rho)}d_jl_j(\rho)$
(see Theorem~\ref{MVlem2}) occurrences. The fourth possibility is
an occurrence of two elements which the first belongs to $C^i$ and
the second belongs to $C^j$ where $i<j$ (Theorem~\ref{MVlem2}
yields every entry of $C^i$ is greater than every entry of $C^j$
for all $i<j$), so there are $\sum_{1\leq i<j\leq f(\rho)}d_id_j$
occurrences. Therefore, by Lemma~\ref{slem1} we have
$$\sign(\pi)=(-1)^{\sum_{j=1}^{f(\rho)}21(C^j)}(-1)^{21(\rho)}(-1)^{\sum_{j=1}^{f(\rho)}d_jl_j(\rho)}(-1)^{\sum_{1\leq i<j\leq f(\rho)}d_id_j},$$
equivalently, $\sign(\pi)=(-1)^{\left(\sum_{1\leq i\leq j\leq
f(\rho)}d_id_j+\sum_{j=1}^{f(\rho)}d_jl_j(\rho)\right)}\cdot
\sign(\rho)\cdot\prod_{j=1}^{f(\rho)}\sign(C^j)$.
\end{proof}
We say the vector $\vv=(v_1,v_2,\ldots,v_n)$ is a {\em binary
vector} if $v_i\in\{0,1\}$ for all $i$, $1\leq i\leq n$. We denote
the set of all binary vectors of length $n$ by $\BB^n$. For any
$\vv\in\BB^n$, we define $|\vv|=v_1+v_2+\cdots+v_n$. For example,
$\BB^2=\{(0,0),(0,1),(1,0),(1,1)\}$ and $|(1,1,0,0,1)|=3$.

Let $\rho$ be any kernel permutations, we denote by $X_a^\rho$
\resp{$Y_a^\rho$} the set of all the binary vectors
$\vv\in\BB^{f(\rho)}$ such that $(-1)^{|\vv|+s(\rho)}=a$
\resp{$(-1)^{|\vv|}=a$}.
For any $\vv\in\BB^{f(\rho)}$, we define
$$z_\rho(\vv)=(-1)^{\sum_{1\leq i<j\leq f(\rho)}v_iv_j+\sum_{j=1}^{f(\rho)}
l_j(\rho)v_j}\sign(\rho).$$ Let $\rho$ be any kernel permutations
and $\vv=(v_1,v_2,\ldots,v_{f(\rho)})\in\BB^{f(\rho)}$, we denote
by $\SS(\rho;\vv)$ the set of all permutations of all sizes with
kernel permutation $\rho$ such that the corresponding cells $C^j$
satisfy $(-1)^{d_j}=(-1)^{v_j}$, in such a context $\vv$ is called
a {\it length argument vector\/} of $\rho$. By definitions, the
following result holds immediately.

\begin{lemma}\label{alem2}
For any kernel permutation $\rho$,
    $$\SS(\rho)=\bigsqcup_{\vv\in\BB^{f(\rho)}}\SS(\rho;\vv).$$
\end{lemma}

Let $\rho$ be any kernel permutations and let
$\vv=(v_1,v_2,\ldots,v_{f(\rho)})$,
$\uu=(u_1,u_2,\ldots,u_{f(\rho)})\in\BB^{f(\rho)}$, we denote by
$\SS(\rho;\vv,\uu)$ the set of all permutations in $\SS(\rho;\vv)$
such that the corresponding cells $C^j$ satisfy $\sign(C^j)=1$ if
and only if $u_j=0$, in such a context $\uu$ is called a {\it
signature argument vector\/} of $\rho$. By Lemma~\ref{alem2}, the
following result holds immediately.
\begin{lemma}\label{alem3}
For any kernel permutation $\rho$,
    $$\SS(\rho)=\bigsqcup_{\vv\in\BB^{f(\rho)}}\SS(\rho;\vv)=\bigsqcup_{\vv\in\BB^{f(\rho)}}
    \bigsqcup_{\uu\in\BB^{f(\rho)}}\SS(\rho;\vv,\uu).$$
\end{lemma}

For any $a,b\in\{0,1\}$ we define
$$H_r(a,b)=\left\{\begin{array}{ll}
                     \frac{1}{2}(E_r(x)+(-1)^aE_r(-x)), &   b=0\\
                     &\\
                     \frac{1}{2}(O_r(x)+(-1)^aO_r(-x)), &   b=1
\end{array}\right..$$
By definitions, the following result holds immediately.
\begin{lemma}\label{alem4}
Let $a,b\in\{0,1\}$. The generating function for all permutations
$\pi$ such that $132(\pi)=r$, $(-1)^{|\pi|}=(-1)^a$, and
$\sign(\pi)=(-1)^b$ is given by $H_r(a,b)$.
\end{lemma}

\section{Main Theorem}
The main result of this note can be formulated as follows. Denote
by $K$ the set of all kernel permutations, and by $K_t$ the set of
all kernel shapes for permutations in $\SS_t$. Let $\rho$ be any
kernel permutation, for any $a,b\in\{0,1\}$ and any
$r_1,\ldots,r_{f(\rho)}$ we define
$$L_{r_1,\ldots,r_{f(\rho)}}^\rho(a,b)=\sum_{\vv\in
X_{(-1)^a}^\rho\;}\sum_{\uu\in Y_{(-1)^b
z_\rho(\vv)}^\rho\;}\prod_{j=1}^{f(\rho)}H_{r_j}(v_j,u_j).$$
%
%
%
%

\begin{theorem}\label{thm}
Let $r\geq1$. For any $a,b\in\{0,1\}$,
\begin{equation}
H_r(a,b)=\sum_{\rho\in
K_{2r+1}}\;\sum_{r_1+\dots+r_{f(\rho)}=r-c(\rho),\, r_j\gs0\;}
L_{r_1,\ldots,r_{f(\rho)}}^\rho(a,b).\label{meq}
\end{equation}
\end{theorem}
\begin{proof}
Let us fix a kernel permutation $\rho\in K_{2r+1}$, a length
argument vector $\vv=(v_1,\ldots,v_{f(\rho)})\in
X_{(-1)^a}(\rho)$, and a signature argument vector
$\uu=(u_1,\ldots,u_{f(\rho)})\in Y_{(-1)^bz_\rho(\vv)}^{f(\rho)}$.
Recall that the kernel $\rho$ of any $\pi$ contains exactly
$c(\rho)$ occurrences of $132$. The remaining $r-c(\rho)$
occurrences of $132$ are distributed between the feasible cells of
the kernel cell decomposition of $\pi$. By Theorem~\ref{MVlem2},
each occurrence of $132$ belongs entirely to one feasible cell,
and the occurrences of $132$ in different cells do not influence
one another.

Let $\pi$ be any permutation such that $132(\pi)=r$,
$\sign(\pi)=(-1)^b$ and $(-1)^{|\pi|}=(-1)^a$ together with a
kernel permutations $\rho$, length argument vector $\vv$, and
signature argument vector $\uu$. Then by Lemma~\ref{alem3}, the
cells $C^j$ satisfy the following conditions:

(1) $v_j=0$ if and only if $d_j$ is an even number,

(2) $u_j=0$ if and only if $\sign(C^j)=1$,

(3) $(-1)^{v_1+\ldots+v_{f(\rho)}+s(\rho)}=(-1)^a$, and

(4) $(-1)^{u_1+\ldots+u_{f(\rho)}}z_{\rho}(\vv)=(-1)^b$.

Therefore, by Lemma~\ref{alem4} this contribution gives
$$x^{s(\rho)}\sum_{r_1+\dots+r_{f(\rho)}=r-c(\rho),\, r_j\gs0}
\;\prod_{j=1}^{f(\rho)}H_{r_j}(v_j,u_j).$$ Hence by
Lemma~\ref{alem3} and \cite[Theorem~1]{MaVa02}, if summing over
all the kernel permutations $\rho\in K_{2r+1}$, length argument
vectors $\vv\in X_{(-1)^a}(\rho)$, and signature argument vectors
$\uu\in Y_{(-1)^bz_\rho(\vv)}^{f(\rho)}$ then we get the desired
result.
\end{proof}
%
Theorem~\ref{thm} provides a finite algorithm for finding $E_r(x)$
and $O_r(x)$ for any given $r\gs0$, since we have to consider all
permutations in $\SS_{2r+1}$, and to perform certain routine
operations with all shapes found so far. Moreover, the amount of
search can be decreased substantially due to the following
proposition.

\begin{proposition}\label{pro1} The only kernel permutation of capacity $r\gs1$
and size $2r+1$ is
$$\rho=2r-1\,2r+1\,2r-3\,2r\,\dots 2r-2j-3\,2r-2j\,\dots 1\,4\,2.$$
Its parameters are given by $s(\rho)=2r+1$, $c(\rho)=r$,
$f(\rho)=r+2$, $\sign(\rho)=-1$, and
$z_\rho(v_1,\ldots,v_{r+2})=(-1)^{\left(1+v_{r+2}+\sum_{1\leq
i<j\leq r+2}v_iv_j\right)}$.
\end{proposition}
\begin{proof}
The first part of the proposition holds by
\cite[Prorposition]{MaVa02}. Besides, by using the form of $\rho$
we get $s(\rho)=2r+1$, $c(\rho)=r$, $f(\rho)=r+2$;
$\sign(\rho)=-1$, and $l_j(\rho)=2r$ for all $j=1,2,\ldots,r+1$
and $l_{r+2}(\rho)=2r+1$. Therefore,
$z_\rho(v_1,\ldots,v_{r+2})=(-1)^{\left(1+v_{r+2}+\sum_{1\leq
i<j\leq r+2}v_iv_j\right)}$.
\end{proof}

By this proposition, it suffices to search only permutations in
$\SS_{2r}$. Below we present several explicit calculations.

\subsection{The case $r=0$} Let us start from the case $r=0$. Observe that
Theorem~\ref{thm} remains valid for $r=0$, provided the left hand
side of Equation~\ref{meq} for $a=b=0$ is replaced by
$H_r(0,0)-1=\frac{1}{2}(E_r(x)+E_r(-x))-1$; subtracting $1$ here
accounts for the empty permutation. So, we begin with finding
kernel shapes for all permutations in $\SS_1$. The only shape
obtained is $\rho_1=1$, and it is easy to see that $s(\rho_1)=1$,
$c(\rho_1)=0$, $f(\rho_1)=2$,
$$X_1(\rho_1)=Y_{-1}=\{(1,0), (0,1)\},\quad X_{-1}(\rho_1)=Y_1=\{(0,0), (1,1)\},$$
and
$$z_{\rho_1}(0,0)=z_{\rho_1}(1,0)=z_{\rho_1}(1,1)=-z_{\rho_1}(0,1)=1.$$
Therefore, Equation~\ref{meq} for $a=b=0$ gives
\begin{equation}\begin{array}{l}
\frac{1}{2}(E_0(x)+E_0(-x))-1=\\
\qquad\qquad=xH_0(1,0)H_0(0,0)+xH_0(1,1)H_0(0,1)+xH_0(1,0)H_0(0,1)+xH_0(1,1)H_0(0,0),
\end{array}\label{ra1}\end{equation}
(Equation~\ref{meq} for $a=1$ and $b=0$ gives
\begin{equation}\begin{array}{l}
\frac{1}{2}(E_0(x)-E_0(-x))=xH_0^2(0,0)+xH_0^2(0,1)+xH_0^2(1,0)+xH_0^2(1,1),
\end{array}\label{ra2}\end{equation}
Equation~\ref{meq} for $a=0$ and $b=1$ gives
\begin{equation}\begin{array}{l}
\frac{1}{2}(O_0(x)+O_0(-x))=\\
\qquad\qquad=xH_0(1,1)H_0(0,0)+xH_0(1,0)H_0(0,1)+xH_0(0,0)H_0(1,0)+xH_0(0,1)H_0(1,1),
\end{array}\label{ra3}\end{equation}
and Equation~\ref{meq} for $a=b=1$ gives
\begin{equation}\begin{array}{l}
\frac{1}{2}(O_0(x)-O_0(-x))=2xH_0(0,1)H_0(0,0)+2xH_0(1,1)H_0(1,0).
\end{array}\label{ra4}\end{equation}
Out present aim to find explicitly $E_0(x)$ and $O_0(x)$, thus we
need the following notation. We define
        $$M_r(x)=E_r(x)-O_r(x)\mbox{ and }F_r(x)=E_r(x)+O_r(x)$$
for all $r\geq0$. Clearly,
$$\begin{array}{l}
H_r(0,0)-H_r(0,1)=\frac{1}{2}(M_r(x)+M_r(-x)),\quad
H_r(0,0)+H_r(0,1)=\frac{1}{2}(F_r(x)+F_r(-x)),\\
H_r(1,0)-H_r(1,1)=\frac{1}{2}(M_r(x)-M_r(-x)),\quad
H_r(1,0)+H_r(1,1)=\frac{1}{2}(F_r(x)-F_r(-x)),\end{array}$$ for
all $r\geq0$. Therefore, by subtracting (respectively; adding)
Equation~\ref{ra3} and Equation~\ref{ra1}, and by subtracting
(respectively; adding) Equation~\ref{ra4} and Equation~\ref{ra2}
we get
$$\left\{\begin{array}{l}
M_0(x)+M_0(-x)=2\\
M_0(x)-M_0(-x)= x(M_0^2(x)+M_0^2(-x))\end{array}\right.\mbox{ and
}\left\{\begin{array}{l}
F_0(x)+F_0(-x)=2+x(F_0^2(x)-F_0^2(-x))\\
F_0(x)-F_0(-x)= x(F_0^2(x)+F_0^2(-x)).\end{array}\right.$$ Hence,
        $$M_0(x)=1+\frac{1-\sqrt{1-4x^2}}{2x}\mbox{ and } F_0(x)=\frac{1-\sqrt{1-4x}}{2x}.$$
\begin{theorem}
{\rm(i)} The generating function for the number of even
permutations avoiding $132$ is given by {\rm(}see {\rm
\cite{SiSc85})}
    $$E_0(x)=\frac{1}{2}\left( \frac{1-\sqrt{1-4x}}{2x}+1+\frac{1-\sqrt{1-4x^2}}{2x} \right).$$
{\rm(ii)} The generating function for the number of odd
permutations avoiding $132$ is given by {\rm(}see {\rm
\cite{SiSc85})}
    $$O_0(x)=\frac{1}{2}\left( \frac{1-\sqrt{1-4x}}{2x}-1-\frac{1-\sqrt{1-4x^2}}{2x} \right).$$
{\rm(iii)} The generating function for the number of permutations
avoiding $132$ is given by {\rm(}see {\rm\cite{Kn})}
    $$F_0(x)=\frac{1-\sqrt{1-4x}}{2x}.$$
\end{theorem}

\subsection{The Case $r=1$}
Since permutations in $\SS_2$ do not exhibit kernel shapes
distinct from $\rho_1$, the only possible new shape is the
exceptional one, $\rho_2=132$. Calculation of the parameters of
$\rho_2$ gives $s(\rho_2)=3$, $c(\rho_2)=1$, $f(\rho_2)=3$,
$$\begin{array}{l}
X_1 (\rho_2)=Y_{-1}=\{(1,0,0), (0,1,0), (0,0,1), (1,1,1)\},\\
X_{-1}(\rho_2)=Y_1=\{(0,0,0), (1,1,0), (1,0,1), (1,1,0)\},
\end{array}$$ and
$$\begin{array}{l}
z_{\rho_2}(0,0,0)=z_{\rho_2}(1,0,0)=z_{\rho_2}(0,1,0)=-z_{\rho_2}(1,1,0)=1,\\
-z_{\rho_2}(0,0,1)=z_{\rho_2}(1,0,1)=z_{\rho_2}(0,1,1)=z_{\rho_2}(1,1,1)=1.
\end{array}$$
Therefore, by Theorem~\ref{thm} we have
$$\left\{\begin{array}{l}
2(H_1(0,0)-H_1(0,1))=\\
\qquad=M_1(x)+M_1(-x)=\frac{x^3}{2}(M_0(-x)-M_0(x))(M_0^2(-x)+M_0^2(x))\\
\\
2(H_1(1,0)-H_1(1,1))=\\
\qquad=M_1(x)-M_1(-x)=2x(M_0(x)M_1(x)+M_0(-x)M_1(-x))\\
\quad\qquad\qquad\qquad\qquad\qquad\qquad\qquad\qquad-\frac{x^3}{2}(M_0(-x)+M_0(x))(M_0^2(-x)+M_0^2(x))),
\end{array}\right.$$
Using the expression for $M_0(x)$ (see the case $r=0$)  we get
        $$M_1(x)=\frac{1}{2}(-1+3x+2x^2)+\frac{1-3x-4x^2+4x^3}{2}(1-4x^2)^{-1/2}.$$
Similarly, if considering the expressions for $H_1(0,0)+H_1(0,1)$
and $H_1(1,0)+H_1(1,1)$ we get
        $$F_1(x)=\frac{1}{2}(x-1)+\frac{1-3x}{2}(1-4x)^{-1/2}.$$
\begin{theorem}
{\rm(i)} The generating function for the number of even
permutations containing $132$ exactly once is given by
    $$E_1(x)=-\frac{1}{2}(1-2x-x^2)+\frac{1-3x}{4}(1-4x)^{-1/2}+\frac{1-3x-4x^2+4x^3}{4}(1-4x^2)^{-1/2}.$$

{\rm(ii)} The generating function for the number of odd
permutations containing $132$ exactly once is given by
    $$O_1(x)=-\frac{1}{2}(x+x^2)+\frac{1-3x}{4}(1-4x)^{-1/2}-\frac{1-3x-4x^2+4x^3}{4}(1-4x^2)^{-1/2}.$$

{\rm(iii)} The generating function for the number of permutations
containing $132$ exactly once is given by {\rm(}see
{\rm\cite{Bo98})}
    $$F_1(x)=\frac{1}{2}(x-1)+\frac{1-3x}{2}(1-4x)^{-1/2}.$$
\end{theorem}

\subsection{The case $r=2$} We have to check the kernel shapes of permutations in $\SS_4$.
Exhaustive search adds four new shapes to the previous list; these
are $1243$, $1342$, $1423$, and $2143$; besides, there is the
exceptional $35142\in\SS_5$. Calculation of the parameters $s$,
$c$, $f$, $z$, $X_a$, $Y_a$ is straightforward, and we get
\begin{theorem}
{\rm(i)} The generating function for the number of even
permutations containing $132$ exactly twice is given by
$$\begin{array}{l}
E_2(x)=\dfrac{1}{2}x(x^3+3x^2-4x-1)+\dfrac{1}{4}(2x^4-4x^3+29x^2-15x+2)(1-4x)^{-3/2}\\
\qquad\qquad\qquad\qquad-\dfrac{1}{4}(16x^7-48x^6-76x^5+64x^4+36x^3-21x^2-5x+2)(1-4x^2)^{-3/2}.
\end{array}$$

{\rm(ii)} The generating function for the number of odd
permutations containing $132$ exactly once is given by
$$\begin{array}{l}
O_2(x)=-\dfrac{1}{2}(x^4+3x^3-5x^2-4x+2)+\dfrac{1}{4}(2x^4-4x^3+29x^2-15x+2)(1-4x)^{-3/2}\\
\qquad\qquad\qquad\qquad+\dfrac{1}{4}(16x^7-48x^6-76x^5+64x^4+36x^3-21x^2-5x+2)(1-4x^2)^{-3/2}.
\end{array}$$

{\rm(iii)} The generating function for the number of permutations
containing $132$ exactly twice is given by {\rm(}see
{\rm\cite{MaVa02})}
    $$F_2(x)=\frac{1}{2}(x^2+3x-2)+\frac{1}{2}(2x^4-4x^3+29x^2-15x+2)(1-4x)^{-3/2}.$$
\end{theorem}

\subsection{The cases $r=3,4,5,6$}
Let $r=3,4,5,6$; exhaustive search in $\SS_6$, $\SS_8$,
$\SS_{10}$, and $\SS_{12}$ reveals $20$, $104$, $503$, and $2576$
new nonexceptional kernel shapes, respectively, and we get

\begin{theorem}
Let $r=3,4,5,6$, then
$$M_r(x)=\frac12\left(A_r(x)+B_r(x)(1-4x^2)^{-r+1/2}\right)\mbox{
and } F_r(x)=\frac12\left(C_r(x)+D_r(x)(1-4x)^{-r+1/2}\right)$$
where
$$\begin{array}{l}\\
A_3(x)=2x^6+10x^5-24x^4-30x^3+23x^2+7x-2,\\ \\
A_4(x)=2x^8+14x^7-46x^6-90x^5+117x^4+85x^3-42x^2-8x+1,\\ \\
A_5(x)=2x^10+18x^9-76x^8-198x^7+360x^6+440x^5-355x^4-171x^3+62x^2+10x-2,\\
\\
A_6(x)=256x^{13}-446x^{12}-618x^{11}+194x^{10}-140x^9+798x^8+1404x^7-1702x^6\\
\qquad\qquad-1430x^5+815x^4+302x^3-88x^2-15x+4,
\end{array}$$

$$\begin{array}{l}\\
B_3(x)=64x^{11}-320x^{10}-800x^9+1216x^8+1124x^7-972x^6-524x^5+312x^4+100x^3-43x^2\\
\qquad\qquad-7x+2,\\ \\
B_4(x)=-256x^{15}+1792x^{14}+6112x^{13}-13120x^{12}-19840x^{11}+22224x^{10}+19054x^9\\
\qquad\qquad-14780x^8-8328x^7+4772x^6+1840x^5-775x^4-197x^3+56x^2+8x-1,\\
\\
B_5(x)=1024x^{19}-9216x^{18}-40064x^{17}+111744x^{16}+228896x^{15}-343264x^{14}-404056x^{13}\\
\qquad\qquad+398712x^{12}+321058x^{11}-234686x^{10}-137468x^9+78480x^8+33896x^7-15400x^6\\
\qquad\qquad-4780x^5+1723x^4 +351x^3-98x^2-10x+2,\\ \\
B_6(x)=524288x^{24}+1175552x^{23}-1593344x^{22}-2324992x^{21}+1162752x^{20}+298112x^{19}\\
\qquad\qquad+2696448x^{18}+4856864x^{17}-7020288x^{16}-7464568x^{15}+6981056x^{14}\\
\qquad\qquad+5445696x^{13}-3868942x^{12}-2335450x^{11}+1324884x^{10}+627306x^9\\
\qquad\qquad-290536x^8-106510x^7+40772x^6+11046x^5-3543x^4-632x^3+176x^2+15x-4,
\end{array}$$
$$\begin{array}{l}\\
C_3(x)=2x^3-5x^2+7x-2,\\ \\
C_4(x)=5x^4-7x^3+2x^2+8x-3,\\ \\
C_5(x)=14x^5-17x^4+x^3-16x^2+14x-2,\\ \\
C_6(x)=42x^6-44x^5+5x^4+4x^3-20x^2+19x-4,
\end{array}$$
and
$$\begin{array}{l}\\
D_3(x)=-22x^6-106x^5+292x^4-302x^3+135x^2-27x+2,\\ \\
D_4(x)=2x^9+218x^8+1074x^7-1754x^6+388x^5+1087x^4,\\ \\
D_5(x)=-50x^{11}-2568x^{10}-10826x^9+16252x^8-12466x^7+16184x^6-16480x^5+9191x^4\\
\qquad\qquad-2893x^3+520x^2-50x+2,\\ \\
D_6(x)=4x^{14}+820x^{13}+32824x^{12}+112328x^{11}-205530x^{10}+141294x^9-30562x^8\\
\qquad\qquad-67602x^7+104256x^6-74090x^5+30839x^4-7902x^3+1230x^2-107x+4.
\end{array}$$

Moreover, for $r=3,4,5,6$,
$$E_r(x)=\frac{1}{4}\left(A_r(x)+C_r(x)+D_r(x)(1-4x)^{-r+1/2}+B_r(x)(1-4x^2)^{-r+1/2}\right)$$
and
$$O_r(x)=\frac{1}{4}\left(A_r(x)-C_r(x)+D_r(x)(1-4x)^{-r+1/2}-B_r(x)(1-4x^2)^{-r+1/2}\right).$$
\end{theorem}
\def\xx{{\bf x}}
\def\yy{{\bf y}}
\section{Further results and open questions}
First of all, let us simplify the expression
$$L_{r_1,\ldots,r_{f(\rho)}}^\rho(a,0)-L_{r_1,\ldots,r_{f(\rho)}}^\rho(a,1),$$
where $a=0,1$, $r_j\gs0$ for all $j$.

\begin{lemma}\label{blema1}
Let $\vv\in\{0,1\}^n$ be any vector, and let $a\in\{1,-1\}$. Then
$$\sum_{\xx\in Y_a}\prod_{j=1}^n H_r(v_j,x_j)-\sum_{\yy\in Y_{-a}}\prod_{j=1}^n H_r(v_j,y_j)
=a\prod_{j=1}^n g_r(j),$$ where
$g_r(j)=H_r(v_j,0)-H_r(v_j,1)=\frac{1}{2}(M_r(x)+(-1)^{v_j}M_r(-x))$
for all $j$.
\end{lemma}
\begin{proof}
Let us define an order on the set $\BB^n$, we say the vector
$\vv<\uu$ if there exists $j$ such that $u_j+v_j=1$, and $u_i=v_i$
for all $i\neq j$. We say the $2^n$ vectors
$\uu^1,\cdots,\uu^{2^n}$ of $\BB^n$ are satisfy the {\em
$\ell$-property} if
$${\bf 0}=(0,0,\ldots,0)=\uu^1<\uu^2<\dots<\uu^{2^n},$$
and we say the $2^n$ vectors $\uu^1,\cdots,\uu^{2^n}$ are satisfy
the {\em c-property} if the vectors
$$(\uu^1_1,\ldots,\uu^1_m),\ldots,(\uu^{2^m}_1,\ldots,\uu^{2^m}_m)$$
are satisfy the $\ell$-property for all $m=1,2,\ldots,n$. For
example, the vectors of $\BB^3$ are satisfy the c-property by
$$(0,0,0)<(1,0,0)<(1,1,0)<(0,1,0)<(0,1,1)<(1,1,1)<(1,0,1)<(0,0,1).$$
First of all, let us prove by induction on $n$ that there exists
an order of the vectors of $\BB^n$ with c-property, for all
$n\geq1$. For $n=1$ the c-property holds with $(0)<(1)$. Suppose
that there exists an order of the vector of $\BB^m$ with the
c-property. Let $\vv^j=(\uu^j_1,\ldots,\uu^j_m,0)$ for all
$j=1,2,\ldots,2^m$ and
$\vv^{2^m+j}=(\uu^{2^m+1-j}_1,\ldots,\uu^{2^m+1-j}_m,1)$ for all
$j=1,2,\ldots,2^m$. By definitions, $\vv^1=(0,0,\ldots,0)$ and
$\vv^1<\cdots<\vv^{2^{m+1}}$, so the $\ell$-property holds for
$m+1$. Hence, by induction on $m$ we get that there exists an
order of the vector of $\BB^n$ with the c-property.

Now we are ready to prove the lemma. Without loss of generality we
can assume that $(0,0,\ldots,0)\in Y_a^n$ (which means $a=1$);
otherwise it is enough to replace $a$ by $-a$. Let
$\xx^1,\ldots,\xx^{2^n}$ all the vectors of $\BB^n$ with the
c-property. Using $(0,0,\ldots,0)\in Y_a^n$ together with the
c-property we get that $\xx^{2i-1}\in Y_a^n$ and $\xx^{2i}\in
Y_{-a}^n$ for all $i=1,2,\ldots,2^{n-1}$. Therefore, for all
$i=1,2,\ldots,2^{n-1}$,
$$\begin{array}{ll}
\sum\limits_{i=1}^{2^{n-1}}\left(\prod\limits_{j=1}^nH_r(v_j,\xx^{2i-1}_j)-\prod\limits_{j=1}^nH_r(v_j,\xx^{2i}_j)\right)
&=\sum\limits_{i=1}^{2^{n-1}}\left((-1)^{i-1}g_r(1)\prod\limits_{j=2}^{n}H_r(v_j,\xx^{2i-1}_j)\right)=\\
&=g_r(1)\sum\limits_{i=1}^{2^{n-2}}\left(\prod\limits_{j=1}^{n-1}H_r(\widetilde{v}_j,\yy^{2i-1}_{j})-\prod\limits_{j=1}^{n-1}H_r(\widetilde{v}_j,\yy^{2i}_{j})\right),
\end{array}$$
where $\yy^p=(\xx^{2p}_2,\ldots,\xx^{2p}_n)$ for all
$p=1,2,\ldots,2^{n-1}$, and $\widetilde{v}=(v_2,v_3,\ldots,v_n)$.
Using the c-property for $\xx^1,\ldots,\xx^{2^n}$ we get that the
vectors $\yy^1,\ldots,\yy^{2^{n-1}}$ are satisfy the c-property in
$\BB^{n-1}$. Hence by induction on $n$ (by definitions, the lemma
holds for $n=1$), we get that the expression equals to
$a\prod_{j=1}^n g_r(j)$.
\end{proof}

As a remark, the vector $(0,0,\ldots,0)\in Y_{z_\rho(\vv)}^{\rho}$
if and only if $z_\rho(\vv)=1$ for any kernel permutation $\rho$
and vector $\vv$. Therefore, by Theorem~\ref{thm} and
Lemma~\ref{blema1} we get the following result.
\begin{theorem}\label{thba1}
Let $a\in\{0,1\}$ and $r\geq0$. Then
$$\begin{array}{l}
\frac{1}{2}(M_r(x)+(-1)^aM_r(-x))-\delta_{r+a,0}=\\
\quad=\sum\limits_{\rho\in
K_{2r+1}}x^{s(\rho)}\sum\limits_{r_1,\ldots,r_{f(\rho)}=r-c(\rho),\,
r_j\gs0}\left(\sum\limits_{\vv\in
X_{(-1)^a}(\rho)}2^{-f(\rho)}z_\rho(\vv)\prod\limits_{j=1}^{f(\rho)}
(M_{r_j}(x)+(-1)^{v_j}M_{r_j}(-x))\right). \end{array}$$
\end{theorem}
As a remark, the above theorem yields two equations (for $a=0$ and
$a=1$) that are linear on $M_r(x)$ and $M_r(-x)$. So,
Theorem~\ref{thba1} provides a finite algorithm for finding
$M_r(x)$ for any given $r\gs0$, since we have to consider all
permutations in $\SS_{2r+1}$, and to perform certain routine
operations with all shapes found so far. Moreover, the amount of
search can be decreased substantially due to the following
proposition which holds immediately by Proposition~\ref{pro1} and
Theorem~\ref{thba1}.

\begin{proposition} Let $r\gs1$, $a\in\{0,1\}$, and
$\rho=2r-1\,2r+1\,2r-3\,2r\,\dots 2r-2j-3\,2r-2j\,\dots 1\,4\,2$.
Then the expression
$$x^{s(\rho)}\sum\limits_{r_1,\ldots,r_{f(\rho)}=r-c(\rho),\,
r_j\gs0}\left(\sum\limits_{\vv\in
X_{(-1)^a}(\rho)}2^{-f(\rho)}z_\rho(\vv)\prod\limits_{j=1}^{f(\rho)}
(M_{r_j}(x)+(-1)^{v_j}M_{r_j}(-x))\right)$$ is given by
$$\sum_{j=a}^{[(r+2)/2]}(-1)^{j-a+1}2^{-r-2}\binom{r+2}{2j+1-a}x^{2r+1}(M_0(x)-M_0(-x))^j(M_0(x)+M_0(-x))^{r+2-j}.$$
\end{proposition}
By this proposition, it is sufficient to search only permutations
in $\SS_{2r}$.
Besides, using Theorem~\ref{thba1} and the case $r=0$ together
with induction on $r$ we get the following result.
\begin{theorem}\label{thba3}
$M_r(x)$ is a rational function on $x$ and $\sqrt{1-4x^2}$ for any
$r\geq0$.
\end{theorem}

In view of our explicit results, we have even a stronger
conjecture.

\begin{conjecture} For any $r\geq1$,
there exist polynomials $A_r(x)$, $B_r(x)$, $C_r(x)$, and $D_r(x)$
with integer coefficients such that
$$\begin{array}{l}
E_r(x)=\frac{1}{4}(A_r(x)+B_r(x))+\frac{1}{4}C_r(x)(1-4x)^{-r+1/2}+\frac{1}{4}D_r(x)(1-4x^2)^{-r+1/2},\\
\\
O_r(x)=\frac{1}{4}(A_r(x)-B_r(x))+\frac{1}{4}C_r(x)(1-4x)^{-r+1/2}-\frac{1}{4}D_r(x)(1-4x^2)^{-r+1/2}.
\end{array}$$
\end{conjecture}


\end{document}